\documentclass[a4paper,11pt]{paper}

\usepackage{amsmath,amsfonts,amssymb,amsthm}
\usepackage[left=25mm,top=30mm,right=25mm,bottom=30mm]{geometry}
\usepackage{courier}
\usepackage{url}

\usepackage[numbers]{natbib}
\newcommand{\norm}[1]{\left\lVert #1 \right\rVert}

\usepackage{fancyhdr}
\theoremstyle{plain}
\newtheorem*{thm}{Theorem}

\theoremstyle{definition}

\newtheorem*{exmp}{Example}

\theoremstyle{remark}

\title{\Large{Comparison of transport map generated by heat flow interpolation and the optimal transport Brenier map}}
\author{\texttt{\large{Anastasiya Tanana}} \thanks{\texttt{E-mail: atanana@utexas.edu}}}

\begin{document}

\maketitle

\begin{abstract}
This note shows that the non-expansive transport map constructed by Y.-H. Kim and E. Milman using heat flow interpolation is in general different from the optimal transport Brenier map.
\end{abstract}

\section{Introduction}

Let $\mu$ and $\nu$ be two Borel probability measures on $\mathbb{R}^n$. A Borel map $T:\mathbb{R}^n\to\mathbb{R}^n$ is said to push $\mu$ forward to $\nu$ (or transport $\mu$ onto $\nu$), denoted by $T_{\#}\mu=\nu$, if $\mu(T^{-1}(\Omega))=\nu(\Omega)$ for every Borel set $\Omega\subset\mathbb{R}^n$, or equivalently, if  for every bounded Borel function $\zeta:\mathbb{R}^n\to\mathbb{R}$ $$\int_{\mathbb{R}^n}\zeta\circ T d\mu=\int_{\mathbb{R}^n}\zeta d\nu.$$

Herein, we consider two pushforward maps: the optimal transport map for quadratic cost function (also known as the Brenier map) and the transport map constructed by Kim and Milman in \citep{kim-milman} using heat flow interpolation. We prove that they are generally different maps, thus answering the question discussed in \citep{kim-milman}.

For this purpose, we consider a Gaussian measure $\mu$ with density
\begin{equation}\label{measure}
\frac{d\mu}{dx}=\frac{\sqrt{\det(A)}}{(2\pi)^{\frac{n}{2}}}\exp\left(-\frac{1}{2}x^\intercal Ax\right),
\end{equation}
where $A$ is a symmetric positive definite matrix, and a Borel probability measure $\nu$ log-concave with respect to $\mu$, that is, $d\nu=\exp(-F)d\mu$ for a convex function $F:\mathbb{R}^n\to\mathbb{R}$. 

The note is organized as follows. In Sections 2 and 3, we recall some facts about the Brenier optimal transport map and sketch the Kim-Milman construction. In Section 4, we show that if we take $\frac{d\nu}{d\mu}=c_0\cdot\exp\left(-\frac{1}{2}x^\intercal Bx\right)$, then we can find $A$ and $B$ such that the two maps do not coincide. These are the probability distributions suggested in Example 6.1 of \citep{kim-milman}, and we show that indeed they can give a counterexample. We also mention a numerical result that suggests that the maps are generally different even in the special case when $\mu$ is the standard normal distribution.

\section{The Brenier map}

The Monge-Kantorovich optimal transport problem with quadratic cost is the problem of finding a minimizer of the functional
$$\int_{\mathbb{R}^n\times\mathbb{R}^n}\norm{x-y}^2d\pi(x,y)$$
over all couplings $\pi$ of $\mu$ and $\nu$, i.e. over all Borel probability measures $\pi$ on $\mathbb{R}^n\times\mathbb{R}^n$ such that for every Borel set $\Omega\subset\mathbb{R}^n$, $\pi(\Omega\times\mathbb{R}^n)=\mu(\Omega)$ and $\pi(\mathbb{R}^n\times \Omega)=\nu(\Omega)$. 

The following result is well-known in optimal transportation theory (for example, see \citep[Theorems 2.12 and 2.32]{villani1}).

\begin{thm}\label{quadratic-cost}
Let $\mu,\nu$ be Borel probability measures on $\mathbb{R}^n$ and assume that $\mu$ is absolutely continuous with respect to the Lebesgue measure. Then there exists a unique, up to a $\mu$-nullset, measurable map $T$ such that $T_{\#}\mu=\nu$  and $T=\nabla\varphi$ for some convex function $\varphi$. If in addition $\mu$ and $\nu$ have finite second order moments, then $(\text{Id}\times \nabla\varphi)_{\#}\mu$ is the unique solution of the Monge-Kantorovich optimal transport problem with quadratic cost.
\end{thm}

The map $\nabla\varphi$, defined up to a $\mu$-nullset, is called the Brenier map.

It was observed by Caffarelli in \citep{caffarelli} that the Brenier map transporting a Gaussian measure $\mu$ onto a probability measure $\nu$ log-concave with respect to $\mu$ is non-expansive (i.e. $1$-Lipschitz).

\section{The Kim-Milman construction}

Kim and Milman's construction produces another non-expansive map transporting log-concave probability measure $\mu$ onto a probability measure $\nu$ log-concave with respect to $\mu$ via semigroup interpolation. Herein, we sketch the construction for the special case of Gaussian $\mu$ defined as in \eqref{measure}.
Consider the second-order differential operator
$$L=\exp\left(\frac{1}{2}x^\intercal Ax\right)\nabla\cdot\left(\exp\left(-\frac{1}{2}x^\intercal Ax\right)\nabla\right)=\Delta-Ax\cdot\nabla.$$
It is known that the solution to
\begin{equation}\label{fokker-planck}
\left\{ \begin{array}{ll}
 \frac{d}{dt}\left(P^A_t(f)\right)=L\left(P^A_t(f)\right)\\
 P^A_0(f)=f
\end{array}
\right.
\end{equation}
(for $f$ smooth and bounded) is given by the Mehler formula (\citep{harge})
$$P^A_t(f)(x)=\int_{\mathbb{R}^n}f\left(\exp(-tA)x+\sqrt{\text{Id}-\exp(-2tA)}y\right)d\mu(y).$$
The family of operators $\left\{P^A_t\right\}_{t\in[0,\infty)}$ defined by \eqref{fokker-planck} is sometimes called the heat semigroup or heat flow with respect to the generator $L$.

Let us now assume that, besides being convex, $F$ is smooth and bounded from below. If we define $d\nu_t=P^A_t(\exp(-F))d\mu$ then $\nu_0=\nu$ and $\nu_t\to\mu$ as $t\to\infty$ in $L^1(\mathbb{R}^n)$. The equation \eqref{fokker-planck} and the definition of $L$ can be used to show that the densities of $\nu_t$ with respect to the Lebesgue measure solve the following transport equation
$$\frac{d}{dt}\left(\frac{d\nu_t}{dx}\right)-\nabla\cdot\left(\left(\frac{d\nu_t}{dx}\right)\nabla\log P^A_t(\exp(-F))\right)=0.$$
It is known for this equation (for example, see \citep[Theorem 5.34]{villani1}) that if there exists a locally Lipschitz family of homeomorphisms $\left\{S_t\right\}_{t\in[0,\infty)}$ solving the initial value problem
\begin{equation}\label{initial_value}
\frac{d}{dt}S_t(x)=w_t(S_t(x)),\quad S_0(x)=x,
\end{equation}
for the velocity field $w_t(x)=-\nabla\log P^A_t(\exp(-F))(x)$, then $S_{t\#}\nu=\nu_t$. It can be shown that if we additionally assume that $F$ is Lipschitz, then such a family of homeomorphisms $\left\{S_t\right\}_{t\in[0,\infty)}$ exists and is unique. Due to smoothess of $w_t$, $S_t$ are in fact diffeomorphisms, and the equation \eqref{initial_value} implies by differentiation that
\begin{equation}\label{jacobian}
\frac{d}{dt}DS_t(x)=Dw_t\big\vert_{S_t(x)}DS_t(x),\quad DS_0\equiv\text{Id}.
\end{equation}
By the Pr\'ekopa-Leindler inequality (see \citep[Theorems 3 and 6]{prekopa}), $-F$ being concave implies $\log P^A_t(\exp(-F))$ is concave and thus $Dw_t=-D^2\log P^A_t(\exp(-F))$ is positive semidefinite at each point. It follows that
$$\frac{d}{dt}(DS_t)^\intercal(x)(DS_t)(x)=(DS_t)^\intercal(x)\left[(Dw_t)^\intercal\big\vert_{S_t(x)}+Dw_t\big\vert_{S_t(x)}\right](DS_t)(x)\ge 0,$$
and therefore $S_t$ are expansions for all $t\ge0$. Their inverses $T_t=S_t^{-1}$ are then non-expansive and can be shown to converge (uniformly on compact sets, up to a subsequence) to a non-expansive map $T$. Since $T_{t\#}\nu_t=\nu$, in the limit $T_{\#}\mu=\nu$. For arbitrary convex $F$, the non-expansive map $T$ transporting  $\mu$ onto $\nu$ is obtained by an approximation argument (see Lemma 3.3 in \citep{kim-milman} and the discussion after it).

\section{Comparison}

In the last section of \citep{kim-milman} Kim and Milman compare their map $T$ with the Brenier map. They give a sufficient condition (6.3) for the two maps to be the same (in particular, when $n=1$, or when $\mu$ and $\nu$ are both radially symmetric, the maps do coincide), but do not manage to show that in general the maps are different. Continuing Example 6.1 in \citep{kim-milman}, we show that there exist Gaussian measures $\mu$ and $\nu$ such that the construction does not give the Brenier map between them.

\begin{exmp}
We consider the special case $\frac{d\mu}{dx}=c\cdot\exp\left(-\frac{1}{2}x^\intercal Ax\right), \frac{d\nu}{d\mu}=c_0\cdot\exp\left(-\frac{1}{2}x^\intercal Bx\right)$, where $A$ and $B$ are symmetric positive definite matrices, and achieve a contradiction assuming that for all such $A$ and $B$ the Kim-Milman map between $\mu$ and $\nu$ coincides with the Brenier map. The matrices $A$ and $B$ giving the contradiction are to be chosen later.

The Mehler formula can be used to obtain that

$$P^A_t\left(c_0\exp\left(-\frac{1}{2}.^\intercal B.\right)\right)(x)=c_t\exp\left(-\frac{1}{2}x^\intercal B_tx\right)$$

for some constants $c_t$ and constant in space symmetric matrices $B_t$ (with $B_0=B$), which are positive semidefinite by the Pr\'ekopa-Leindler inequality and decay exponentially to $0$ as $t\to\infty$. We obtain $w_t(x)=-\nabla\log P^A_t(\exp(-\frac{1}{2}.^\intercal B.))(x)=B_tx$ and $Dw_t\equiv B_t$. For such matrices $B_t$, Picard-Lindel\"of-type argument and integral Gronwall's lemma imply that both ordinary differential equations \eqref{initial_value} and \eqref{jacobian} have unique solutions well-defined for each $x\in\mathbb{R}^n$ and for all $t\in[0,\infty)$. Clearly, $S_t$ are then linear maps given by multiplication by constant in space matrices $DS_t$.

The explicit expression for $\nu_t$ is

$$d\nu_t=d_t\exp\left(-\frac{1}{2}x^\intercal (A+B_t)x\right)dx,$$

where $d_t=\frac{\sqrt{\det(A+B_t)}}{(2\pi)^{n/2}}$ are the normalizing constants. Hence, $\nu_t$ are also Gaussian and log-concave with respect to $\mu$. Fix $t\ge0$ and consider Kim and Milman's construction for measures $\mu$ and $\tilde{\nu}=\nu_t$. Notice that the flow of measures interpolating between $\tilde{\nu}$ and $\mu$ is the time-shifted initial flow $\nu_t$:

$$d\tilde{\nu}_s=P^A_s\left(P^A_t\left(c_0\exp\left(-\frac{1}{2}.^\intercal B.\right)\right)\right)d\mu=P^A_{s+t}\left(c_0\exp\left(-\frac{1}{2}.^\intercal B.\right)\right)d\mu=d\nu_{t+s},\quad\forall s\ge 0.$$

This is a consequence of the semigroup property for $P^A$: $P^A_s\circ P^A_t=P^A_{s+t}$ for all $s,t\ge 0$, which can be derived, for example, from the Mehler formula.
For the same reason, the corresponding velocity field $\tilde{w}_s=-\nabla\log P^A_s(P^A_t(c_0\exp(-\frac{1}{2}.^\intercal B.)))$ is the time-shifted initial velocity field: $\tilde{w}_s=w_{t+s}$. This implies that the flow of diffeomorphisms $S_s$ along $w_s$ and the flow of diffeomorphisms $\tilde{S}_s$ along $\tilde{w}_s$ ($\tilde{S}_{s\#}\tilde{\nu}=\tilde{\nu}_s$) satisfy

$$S_{t+s}=\tilde{S}_s\circ S_t,\quad\forall s\ge 0.$$

Then the inverse diffeomorphisms $T_{s}=S_s^{-1}$ and $\tilde{T}_s=\tilde{S}_s^{-1}$ satisfy the relation

\begin{equation}\label{T}
\tilde{T}_s=S_t\circ T_{t+s},\quad\forall s\ge 0.
\end{equation}

Denote by $T_{0,opt}$ the Brenier map between $\mu$ and $\nu$, and by $T_{t,opt}$ the Brenier map between $\mu$ and $\tilde{\nu}=\nu_t$. By our assumption, $T_{t+s}\to T_{0,opt}$ and $\tilde{T}_s\to T_{t,opt}$ as $s\to\infty$. In particular, taking the limit as $s\to\infty$ in \eqref{T} gives
\begin{equation}\label{formula}
T_{t,opt}=S_t\circ T_{0,opt},\quad\forall t\ge0.
\end{equation}

Since $\nu_t$ and $\mu$ are Gaussian, the Brenier map between $\nu_t$ and $\mu$ is given explicitly (e.g. \citep[Example 1.7]{mccann2}) by multiplication by the symmetric positive definite matrix
$$A^{1/2}(A^{1/2}(A+B_t)A^{1/2})^{-1/2}A^{1/2}.$$

Therefore, the Brenier map $T_{t,opt}$ between $\mu$ and $\nu_t$, being the unique map pushing $\mu$ forward to $\nu_t$ which is a gradient of a convex function, should be given by multiplication by the inverse of this matrix, i.e.
$$DT_{t,opt}(x)=A^{-1/2}(A^{1/2}(A+B_t)A^{1/2})^{1/2}A^{-1/2},\quad\forall x\in\mathbb{R}^n.$$
Recall that \eqref{jacobian} becomes the following matrix differential equation (identical for all $x$):
$$\frac{d}{dt}DS_t=B_t(DS_t),\quad DS_0=\text{Id}.$$
Multiplying this ODE from the right by the matrix $DT_{0,opt}$, we obtain from \eqref{formula} that $DT_{t,opt}$ satisfy the ODE $\frac{d}{dt}DT_{t,opt}=B_t(DT_{t,opt})$ as well. In particular, since $DT_{t,opt}$ are symmetric, $B_t(DT_{t,opt})$ should be symmetric for all $t$. Consider $t=0$: 
$$B_0(DT_{0,opt})=BA^{-1/2}(A^{1/2}(A+B)A^{1/2})^{1/2}A^{-1/2}.$$
This matrix is symmetric if and only if $C=A^{1/2}BA^{-1/2}(A^{1/2}(A+B)A^{1/2})^{1/2}$ is symmetric. But it is easy to find matrices $A$ and $B$ such that $C$ is not symmetric. For example, take
$$A=\begin{pmatrix}4&0\\0&1\end{pmatrix},
\quad B=\begin{pmatrix}2&1\\1&3\end{pmatrix}.$$
In this case
$$A^{1/2}BA^{-1/2}=\begin{pmatrix}2&2\\0.5&3\end{pmatrix},\quad
A^{1/2}(A+B)A^{1/2}=\begin{pmatrix}24&2\\2&4\end{pmatrix},$$
and one can compute that
$$C\approx\begin{pmatrix}
10.4 & 4.5\\
3.3 & 6.1
\end{pmatrix}.$$
\qed
\end{exmp}

\textbf{The case of standard normal distribution $\mu$}

When $\mu$ is the standard normal distribution ($A=\text{Id}$) and $\nu$ is Gaussian, the Kim-Milman construction does give the Brenier map (\citep[Section 6]{kim-milman}). However, our numerical result suggests that this does not hold for general $d\nu=\exp(-F)d\mu$ with convex function $F$. We consider the case $n=2$, $F(x)=F(x_1,x_2)=x_1^4+x_2^4+(x_1+x_2)^2$ and the starting point $x=(-0.5,0)$. Our numerical solution of \eqref{initial_value} yielded $S_\infty(x)\approx(-1.054,-0.231)$, while the numerical solution to \eqref{jacobian} converged as $t\to\infty$ to a non-symmetric matrix approximately equal to
$$\begin{pmatrix}
2.303&0.441\\
0.467&2.013
\end{pmatrix},$$
meaning that $S_\infty$ is unlikely to be a gradient of a convex function at point $x$. To obtain this approximations, we used the explicit Euler method for both ODE's with terminal time $T=30$ and the following time step sizes:
\begin{table}[h!]
\centering
 \begin{tabular}{|| c c ||} 
 \hline
 time interval &  $\Delta t$ \\ [0.5ex] 
 \hline\hline
[0,0.1] & 0.00002 \\ 
\hline
[0.1,0.5] & 0.00005 \\
\hline
[0.5,1] & 0.0002 \\
\hline
[1,3] & 0.0005 \\
\hline
[1,5] & 0.002 \\
\hline
[5,30] & 0.005 \\
 \hline
 \end{tabular}
\end{table}

To approximate the semigroup $P^A_t$, which becomes the Ornstein-Uhlenbeck semigroup when $A=\text{Id}$, \texttt{scipy.integrate.nquad} was used.
\vspace{12pt}

\textbf{Acknowledgement.} I thank my Master's thesis advisor Joe Neeman for numerous helpful discussions of the topic.

\vspace{12pt}



\footnotesize

\begin{thebibliography}{1}
\bibitem{caffarelli}
	L. A. Caffarelli. Monotonicity properties of optimal transportation and the {FKG} and related inequalities.
 	\textit{Communications in Mathematical Physics}, 214(3):547-563, 2000.
	
\bibitem{harge}
	G. Harg\'e. A particular case of correlation inequality for the {G}aussian measure.
	\textit{The Annals of Probability}, 27(4), 1999.

\bibitem{kim-milman}
	Y.-H. Kim and E. Milman. A generalization of {C}affarelli's contraction theorem via (reverse) heat flow.
	\textit{Mathematische Annalen}, 354(3):827-862, 2012.

\bibitem{mccann2}
	R. J. McCann. A convexity principle for interacting gases. 
	\textit{Advances in Mathematics}, 128:153-179, 1997.
	
\bibitem{prekopa}
	A. Pr\'ekopa. On logarithmic concave measures and functions.
	\textit{Acta Scientiarum Mathematicarum}, 34:335-343, 1973.

\bibitem{villani1}
	C. Villani. \textit{Topics in optimal transportation}, volume 58 of \textit{Graduate Studies in Mathematics}.
	American Mathematical Society, 2003.
\end{thebibliography}
\end{document}